\documentclass{amsart}

\usepackage{amsmath,amssymb,amsthm,graphicx,slashed,subcaption}%
\usepackage{amsmath}%
\usepackage{amsfonts}%
\usepackage{amssymb}%
\usepackage{float}
\usepackage{tikz,pgfplots}
\pgfplotsset{compat=1.11}

\usetikzlibrary{patterns}
\usepackage{mathpple}

\usepackage{hyperref,enumerate}

\usepackage{colortbl}

\usepackage[all,cmtip]{xy}

\definecolor{Red}{rgb}{0.6,0,0}

\usepackage{graphicx}
\providecommand{\U}[1]{\protect\rule{.1in}{.1in}}
\newtheorem{thm}{Theorem}[]
\newtheorem{corl}[thm]{Corollary}

\newtheorem{lma}[thm]{Lemma}

\newtheorem{prop}[thm]{Proposition}
\newtheorem{defn}[thm]{Definition}
\newtheorem{ex}[thm]{Example}
\newtheorem{rem}[thm]{Remark}
\def\tilde{\widetilde}

\def\C{\mathbb{C}}

\DeclareMathOperator{\CCl}{\mathbb C l}

\DeclareMathOperator{\dom}{Dom}
\def\E{\mathcal{E}}
\def\env{\mathrm{env}}
\def\epsilon{\varepsilon}

\def\H{\mathcal{H}}

\DeclareMathOperator{\ind}{Ind}

\def\id{\textup{id}}
\def\K{\mathcal{K}}
\def\L{\mathcal{L}}

\def\min{\mathrm{min}}

\def\R{\mathbb{R}}

\DeclareMathOperator{\Sig}{Sig}

\newcommand{\Toep}[1]
           {\textup{Toep}_{#1}(\C)}
\def\U{\mathcal{U}}
\def\V{\mathcal{V}}
\def\Z{\mathbb{Z}}

\newtheoremstyle{commentstyle}
  {0.2cm}{0.2cm}
  {\sf}
  {0cm}
  {\bfseries}{ }
  {0cm}
  {\thmname{#1}\thmnumber{ #2}:\thmnote{ #3}}

\theoremstyle{commentstyle}
\newtheorem{mycomment}{Comment}

\title{Higher K-groups for operator systems}

\author{Walter D. van Suijlekom}

\address{Institute for Mathematics, Astrophysics and Particle Physics, Radboud
University Nijmegen, Heyendaalseweg 135, 6525 AJ Nijmegen, The Netherlands.
}
\email{waltervs@math.ru.nl}
\date{November 5, 2024}

\begin{document}

\begin{abstract}
  We extend our previous definition of K-theoretic invariants for operator systems based on hermitian forms to higher K-theoretical invariants. We realize the need for a positive parameter $\delta$ as a measure for the spectral gap of the representatives for the K-theory classes. For each $\delta$ and integer $p \geq 0$ this gives operator system invariants $\V_p^\delta(-,n)$, indexed by the corresponding matrix size.

  The corresponding direct system of these invariants has a direct limit that possesses a semigroup structure, and we define the $K_p^\delta$-groups as the corresponding Grothendieck groups. This is an invariant of unital operator systems, and, more generally, an invariant up to Morita equivalence of operator systems. Moreover, there is a formal periodicity that reduces all these groups to either $K_0^\delta$ or $K_1^\delta$. We illustrate our invariants by means of the spectral localizer.  
  \end{abstract}

\maketitle

\tableofcontents

\section{Introduction}
We continue our quest for a definition of K-theory for operator systems. In a previous work, we have defined a $K_0$-group which we proved to be an invariant of operator systems up to complete order isomorphism. Moreover, we showed that it is stably equivalent, which by \cite{EKT21} amounts to invariance under Morita equivalence of unital operator systems. 

In this work, we extend this and define higher K-groups. At the same time, with the applications to the spectral localizer of \cite{LS18b} in mind, we realize the need for a parameter $\delta$ that indicates (and quantifies the size of) a spectral gap in the elements that represent K-theory classes. This also solves the problem one is facing when considering the space of all invertible elements in a finite-dimensional operator system, namely, that this space is contractible.

The definition of $K_0^\delta$ extends the group $K_0$ considered in \cite{Sui24} which is obtained as a special case by setting $\delta =0$. The groups $K_1^\delta(E)$ are defined, and shown to be stably equivalent as well. Moreover, we may introduce higher K-groups $K_p^\delta(E)$ by exploiting graded Clifford algebras, for which we show a formal periodicity to hold, thus reducing them to $K_0^\delta$ and $K_1^\delta$ for $p$ even and odd, respectively. 

As an illustrative example, we obtain the spectral localizer as an index pairing $\ind_D^\delta: K_p^\delta(E) \to \Z$ in the even $(p=0)$ and odd $(p=1)$ case, respectively.

\subsection*{Acknowledgements}
I would like to thank Alain Connes and Kristin Courtney for useful discussions. I am grateful to Yuezhao Li for his suggestions and comments.

\section{Background on operator systems}

We start by briefly recalling the theory of operator systems, referring to \cite{ER00,Pau02,Pis03,Bla06} for more details.

A {\em unital operator system} $(E,e)$ is a matrix-ordered $*$-vector space $E$, equipped with an Archimedean order unit $e$. A map $\phi: E \to F$ between operator systems determines a family of maps $\phi^{(n)} : M_n(E) \to M_n(F)$ given by $\phi^{(n)}([x_{ij}]) = [\phi(x_{ij})]$. A map $\phi: E \to F$ between unital operator systems is called {\em completely positive} if each $\phi^{(n)}$ is positive ($n \geq 1$). We also abbreviate  completely positive by {\em cp}, and unital completely positive by {\em ucp}. 

A {\em dilation} of a ucp map $\phi: E \to B(\H)$ of a unital operator system is a ucp map $\psi: E \to B(\K)$, where $\K$ is a Hilbert space containing $\H$ such that $P_\H \psi(x)|_\H = \phi(x)$ for all $x \in E$. 
The ucp map $\phi$ is called {\em maximal} if every dilation of $\phi$ is obtained by attaching a direct summand.

A non-zero cp map $\phi:E \to B(\H)$ is said to be {\em pure} if the only cp maps satisfying $0 \leq \psi \leq \phi$ are scalar multiples of $\phi$.

We may view $E$ as a concrete operator system in the $C^*$-algebra $C^*(E)$ it generates; in this case, we say that a ucp map $\phi: E \to B(\H)$ has the {\em unique extension property} if it has a unique ucp extension to $C^*(E)$ which is a $*$-representation. If, in addition, the $*$-representation is irreducible, it is called a {\em boundary representation} \cite{Arv69}. The following result is well-known in the literature \cite{Arv69,MS98,Arv08,Kle14,DK15} (see also \cite{Sui24}). 
\begin{prop}
  \label{prop:bdry-pure-max}
  Let $\phi: E \to B(\H)$ be a ucp map. Then
  \begin{enumerate}
  \item $\phi $ is maximal if and only if it has the unique extension property.
    \item $\phi$ is pure and maximal if and only if it is a boundary representation. 
    \end{enumerate}
\end{prop}

\section{$\delta$-gapped K-theory for unital operator systems: $K_0^\delta$ and $K_1^\delta$}
\label{sect:K1}

Even though in an operator system we cannot speak about invertible  elements, we may use the pure and maximal ucp maps to quantify the extent to which an element is (non)singular. 
\begin{defn}
\label{defn:nondeg}
Let $(E,e)$ be a unital operator system and let $\delta > 0$. An element $x \in M_n(E)$ is called $\delta$-{\em singular} if for all $s \in (0,\delta)$ there exists a $g>0$ such that for all pure and maximal ucp maps  $\phi: E \to B(\H)$ we have
\begin{equation}
\left  | \phi^{(2n)} \left( s \cdot e_{2n} +  \begin{pmatrix}0  & x \\ x^* &  0  \end{pmatrix} \right)\right|
    \geq g \cdot \id_\H^{\oplus 2 n}
    \label{eq:nondeg}
    \end{equation}
The largest real number $g>0$ such that \eqref{eq:nondeg} holds is called the {\em $s$-gap} of $x$. 
\end{defn}
We will write $G^\delta(E,n)$ for all $\delta$-singular elements in $M_n(E)$, and $H^\delta(E,n)$ for the $\delta$-singular self-adjoint  elements in $M_n(E)$. We extend this definition of $\delta$-singular to the limiting case $\delta =0 $ to be those elements for which \eqref{eq:nondeg} holds with $s =0$; in the self-adjoint case these are precisely the non-singular hermitian forms considered previously in \cite{Sui24}.

\begin{prop}
    An element $x \in M_n(E)$ is $\delta$-singular if and only if 
$$\begin{pmatrix} s  & \imath_E^{(n)}(x) \\ \imath_E^{(n)}(x^*) & s
  \end{pmatrix}
  $$
  are invertible elements in $M_{2n}(C^*_\env(E))$ for all $s \in (0,\delta)$, where $\imath_E : E \to C^*_\env(E)$ is the $C^*$-envelope. 
  \end{prop}
\proof
In \cite{DK15} the $C^*$-envelope of $E$ is constructed as the direct sum of all boundary representations $(\H_\sigma, \sigma)$:
$$
\imath_E : E \to \bigoplus_\sigma B( \H_\sigma ).
$$
Now if $x \in M_n(E)$ then a matrix such as the above is invertible if and only if its image $\imath_E^{(n)}(x)$ is bounded from below by $g \cdot \id_{\oplus_\sigma (\H_\sigma)}$. But this holds if and only if \eqref{eq:nondeg} holds for all boundary representations $\sigma$. Since by Proposition \ref{prop:bdry-pure-max} a ucp map is a boundary representation if and only if it is pure and maximal, the result follows.
\endproof

We can now translate $\delta$-singularity into spectral properties of $\imath_E^{(n)}(x)$ in $C^*_\env(E)$. For this consider the following spectrum
$$
\Sigma_x := \sigma  \begin{pmatrix} 0  & \imath_E^{(n)}(x) \\ \imath_E^{(n)}(x^*) & 0 \end{pmatrix} \equiv \left\{ \lambda \in \C \mid \left(\begin{smallmatrix} - \lambda  & \imath_E^{(n)}(x) \\ \imath_E^{(n)}(x^*) & -\lambda \end{smallmatrix} \right) \text{ is not invertible } \right\} 
$$
This subset of $\R$ is symmetric around $0$ since we may unitarily conjugate
$$
\begin{pmatrix} 1 & 0 \\ 0 & -1 \end{pmatrix}
\begin{pmatrix} -\lambda  & \imath_E^{(n)}(x) \\ \imath_E^{(n)}(x^*) & -\lambda \end{pmatrix}
\begin{pmatrix} 1 & 0 \\ 0 & -1 \end{pmatrix} = -  \begin{pmatrix} \lambda  & \imath_E^{(n)}(x) \\ \imath_E^{(n)}(x^*) & \lambda \end{pmatrix}
$$
The notion of $\delta$-singularity can now be captured by the property that $\pm s \notin \Sigma_x$ for all $s \in (0,\delta)$, or, equivalently: 

\begin{corl}
  Let $\delta>0$.
  An element $x \in M_n(E)$ is $\delta$-singular iff 
  $\Sigma_x \subseteq (-\infty, -\delta] \cup \{ 0 \} \cup [\delta,\infty)$.
      \label{corl:gap}
\end{corl}
Our notion of $\delta$-singularity thus allows to detect the presence of a so-called {\em spectral gap} of size $\delta$. This terminology is borrowed from physics where the spectral gap is the energy difference between the ground state and the first excited state. For this reason we also refer to such elements as being {\em $\delta$-gapped}. Moreover, we can find an explicit lower bound for the $s$-gap, since
\begin{equation}
\left| \begin{pmatrix} s  & \imath_E^{(n)}(x) \\ \imath_E^{(n)}(x^*) & s
\end{pmatrix} \right|\geq \min \{ \lambda : \lambda \in s+ \Sigma_x \} = \min\{s,\delta-s\}
\label{eq:s-gap}
\end{equation}

In the self-adjoint case, we have the following result. 
\begin{corl}
  Let $\delta>0$.   A self-adjoint element $x \in M_n(E)$ is $\delta$-singular iff $\imath_E^{(n)}(x)$ has spectral gap $\delta$, {\em i.e.} $  \sigma(\imath_E^{(n)}(x)) \subseteq  (-\infty, -\delta] \cup \{ 0 \} \cup [\delta,\infty)$.
\end{corl}
\proof
This follows by a similarity transformation:
\begin{equation*}
\begin{pmatrix} s  & \imath_E^{(n)}(x) \\ \imath_E^{(n)}(x) & s \end{pmatrix}
\sim \begin{pmatrix} s + \imath_E^{(n)}(x)  &  0  \\ 0 & s-  \imath_E^{(n)}(x) \end{pmatrix}   .
\end{equation*}
Indeed, the matrix on the right-hand side is invertible iff $ \pm s \notin \sigma(\imath_E^{(n)}(x))$ for all $s \in (0,\delta)$.
\endproof



\subsection{The $K_0^\delta$-groups}
The $K_0$-group of non-singular self-adjoint elements has been treated in detail in \cite{Sui24}; let us consider here the changes when $\delta$-singular, or $\delta$-gapped self-adjoint elements are considered.

\begin{defn}
  \label{defn:class-herm}
  Let $\delta \geq 0$ and $x,x' \in H^\delta(E,n)$. We say that $x \sim_n x' $ if there exists $\tilde x \in H^\delta(C([0,1]) \otimes E),n)$ such that
  $$
\tilde x(0) = x ; \qquad \tilde x(1) = x'
  $$
We denote the equivalence class of $x \in H^\delta(E,n)$ by $[x]_n$, and the set of all such equivalence classes in $H^\delta(E,n)$ by $\V_0^\delta(E,n)$, or, equivalently,
$$
\V_0^\delta(E,n) =H^\delta(E,n)/_{\sim_n}.
$$
\end{defn}
All these sets $\V_0^\delta(E,n)$ for $n \geq 1$ are invariants of unital operator systems: 
\begin{prop}
  If $E$ and $F$ are completely order isomorphic then $\V_0^\delta(E,n) \cong \V_0^\delta(F,n)$, for any $n \geq 1$.
\end{prop}
\proof
This follows directly from the fact that in this case $C^*_\env(E) \cong C^*_\env(F)$ via a unital $*$-isomorphism. 
\endproof

If we combine this structure at each level $n$ with the direct sum of $\delta$-gapped self-adjoint elements, we obtain a semigroup structure as follows. Consider the direct system of sets $(\V_0^\delta(E,n), \imath_{nm})$ where for $m \geq n$
\begin{align}
  \label{eq:dir-syst}
  \imath_{nm}: \V_0^\delta(E,n) &\to \V_0^\delta(E,m)\\
[x]_n & \mapsto [x \oplus e_{m-n}]_m, \nonumber
\end{align}
We denote the direct limit of the direct system \eqref{eq:dir-syst} by $\varinjlim \V_0^\delta(E,n)$. A more explicit description is given as follows: for $x \in H^\delta(E,n)$ and $x' \in H^\delta(E,n')$ we write $x \sim x'$ if there exists a $k \geq n,n'$ such that $ x \oplus e_{k-n} \sim_k x'\oplus e_{k-n'} $ in $H^\delta(E,k)$. We will write $[x]_E$, or simply $[x]$, for the equivalence class corresponding to $x \in H^\delta(E,n)$ and 
$\V_0^\delta(E) := \amalg_n \V_0^\delta(E,n)/_\sim$ for the corresponding set of equivalence classes. The following is then clear from the definition of the direct limit. 
\begin{prop}
  The set $\V_0^\delta(E)$ is the direct limit $\varinjlim \V_0^\delta(E,n)$ of the direct system \eqref{eq:dir-syst}. Moreover, it is a semigroup when equipped with the direct sum $[x]+ [x'] = [x \oplus x']$ and identity element $0 = [e]$. 
\end{prop}
The $K_0^\delta$-group $K_0^\delta(E)$ for a unital operator system $E$  is now defined to be the Grothendieck group of $\V_0^\delta(E)$; we will refer to it as the $\delta$-gapped $K_0$-group of $E$. In \cite{Sui24} we have shown that for $\delta=0$, that is, for non-singular self-adjoint elements, this group coincides with the $C^*$-algebraic $K_0$-group in the case of a $C^*$-algebra, that it behaves well with respect to direct sums, and that it is invariant under Morita equivalence of unital operator systems. For later use, we record the following result also in the case of general $\delta \geq 0$. The proof follows {\em mutatis mutandis} from that of \cite[Theorem 18]{Sui24}
\begin{prop}[Stability of $K_0^\delta$]
  \label{prop:stab}
  Let $E$ be a unital operator system, $\delta >  0$, and let $N$ be a natural number. Then $\V^\delta_0(E)$ is isomorphic to $\V^\delta_0(M_N(E))$ (and so are the corresponding $K_0^\delta$-groups).
\end{prop}

\subsection{$K_1^\delta$ for unital operator systems}
For the analogue of $K_1$ for operator systems we drop the assumption of being self-adjoint, and consider $\delta$-singularity in the sense of Definition \ref{defn:nondeg}, or, in view of Corollary \ref{corl:gap}, the presence of a spectral gap of size $\delta$. Let us start with some motivating examples.

\begin{ex}
  \label{ex:delta-nonsing}
  \begin{enumerate}[(i)]
\item (Almost unitaries and quantitative $K$-theory) 
  In \cite{OY15} the quantitative $K_1$-group  of a filtered $C^*$-algebra $A = (A_r)$ was defined  in terms of ${\epsilon}$-$r$-unitaries, {\em i.e.} elements $u \in M_n(A_r)$ such that $\| u^* u  - 1\|, \| u u^*  - 1\| < {\epsilon}$. This implies that 
  $$
\begin{pmatrix} 0 & u \\ u^* & 0 \end{pmatrix}^2 > (1-{\epsilon}) .
$$
In other words, the spectrum of the matrix on the left-hand side has spectrum contained in $(-\infty, -\delta] \cup [\delta , \infty)$ for  all $\delta < (1-{\epsilon})^{1/2}$. We conclude that $u \in G^\delta(A_r,n)$ for all those $\delta$ (including $\delta=0$).

\item (Unitaries in operator spaces and systems) In \cite{BN11} unitaries in operator spaces and operator systems where characterized. These turn out to be unitaries in matrices with entries in the $C^*$-envelope, so that we find that they define elements in $G^\delta(E,n)$ as long as $0 \leq \delta <1$.

\item (The odd spectral localizer)
  The odd spectral localizer \cite{LS18b} is defined in terms of a spectral compression $x = P T P$ of an invertible element $T$ in a $C^*$-algebra $A \subseteq B(\H)$ by a projection in $\H$. We then have 
 \begin{align}
   \begin{pmatrix}0  & x \\ x^* & 0  \end{pmatrix}^2 =P  \begin{pmatrix} T T^* \\ 0 & T^*T   \end{pmatrix}^2 P + P  \begin{pmatrix}[P, T][P, T^*] \\ 0 & [P,T^*][P,T]  \end{pmatrix}^2 P
  \label{eq:trick}
 \end{align}
 Since $T$ is invertible, it is bounded below by some $g>0$; the same applies to $T^*$. Consequently, we find from \eqref{eq:trick}
 $$
\begin{pmatrix} 0  & x \\ x^* & 0  \end{pmatrix}^2 \geq (g^2 - \eta^2)P
$$
in terms of $\eta = \| [P,T]\|$. So, provided $\delta <(g^2-\eta^2)^{1/2}$ and $\eta < g$  we find that $x \in G^\delta(PAP,1)$.

\item (Toeplitz matrices)
  \label{ex:toeplitz}
  Consider the unitary $u$ given by the bilateral shift operator $\ell^2(\Z)$. If we compress this with a projection onto $\text{span}_\C \{ e_1, \ldots e_n\}$ in terms of the canonical orthonormal basis of $\ell^2(\Z)$, we obtain the following Toeplitz matrix:
  $$
  PuP = \begin{pmatrix} 0 & 1 & \cdots & \cdots & 0 \\ 0 & 0 & 1 & \cdots & 0  \\
    \vdots & \vdots & \ddots & \ddots & 0 \\ 0 & 0 & \cdots & 0& 1\\ 0 & 0 & \cdots & 0 & 0 \end{pmatrix}
  $$
  This is a candidate for index pairings in the spectral localizer (for the circle, \cite{LS18b,Hek21}), and as such we expect it to represent an element in the Toeplitz operator system. Note that even though $u$ is invertible, $PuP$ is not invertible (in fact, for this case we have $g = \eta = 1$ in the previous example). However, the eigenvalues of
  $$
\begin{pmatrix} s   & PuP \\ Pu^*P & s   \end{pmatrix} 
  $$
are shifted from those of $PuP$ to become equal to $-1+s, s, 1+s$. Hence, $PuP$ is $\delta$-singular provided $0< \delta <1$.  
\end{enumerate}  
\end{ex}

\begin{defn}
  \label{defn:class-Gn}
  Let $ \delta \geq 0$ and $x,x' \in G^\delta(E,n)$. We say that $x \sim_n x' $ if there exists an element $\tilde x \in G^\delta(C([0,1]) \otimes E),n)$ such that
  $$
\tilde x(0) = x ; \qquad \tilde x(1) = x'
  $$
We will write $\V_1^\delta(E,n) = G^\delta(E,n) /_{ \sim_n}$ for the set of equivalence classes of elements in $G^\delta(E,n)$. 
\end{defn}

\begin{rem}
  \label{rem:steffen}
  Let us come back to the parameter $\delta$ in the definition of $\delta$-singularity (Definition \ref{defn:nondeg}). In fact, we could have set $\delta=0$ from the start (as we did with hermitian forms in \cite{Sui24}) and have simply considered all invertible elements in $\imath_E(E) \subseteq C^*_\env(E)$ up to homotopy equivalence. However, in most cases of interest ---including the applications to the spectral localizer ({\em cf.} Example \ref{ex:delta-nonsing}.\ref{ex:toeplitz} and Section \ref{sect:spectr-loc} below)--- such an invariant becomes trivial. Indeed, if $E$ has  a finite-dimensional $C^*$-envelope $M_d(\C)$ then we claim that $\imath_E(E) \cap M_d(\C)^\times$ is contractible to the point $\imath_E(e) = \mathbb{I}_d$ (and the same applies to $M_n(E)$). Namely, for any $x  \in \imath_E(E) \cap M_d(\C)^\times$ choose $z \in S^1$ so that no eigenvalue of $x$ lies on the ray through the origin and $z$ (since $x$ has finitely many eigenvalues, this is always possible).  For any $t \in [0,1]$ we define  $\gamma(t) = t x + z(1-t) \mathbb I_d$. Then $\det \gamma(t) = 0$ iff $z(1-1/t)$ is an eigenvalue of $x$. By our choice of $z$ we find that $\det \gamma(t) \neq 0$ for all $t \in [0,1]$. Moreover, since complex linear combinations of $x$ and $ \mathbb I_d$ are contained in $\imath_E(E)$ it follows that the path $\gamma(t)$ is contained in $\imath_E(E) \cap M_d(\C)^\times$. More generally, we may conclude that $\imath_E^{(n)} (M_n(E)) \cap M_{nd}(\C)^\times$ is contractible to a point for any $n \geq 1$. 
  \end{rem}

Again, we have the following invariance properties of the sets $\V_1^\delta(E,n)$:
\begin{prop}
  If $E$ and $F$ are unital operator systems. If $\phi:E \to F$ is a ucp map for which there exists a $*$-homomorphism $\tilde\phi$ that makes the following diagram commute,
  \begin{equation}
  \xymatrix {
    E \ar[r]^{\phi} \ar[d]_{\imath_E} & F  \ar[d]_{\imath_F}\\
    C^*_\env(E) \ar[r]^{\tilde \phi} & C^*_\env(F)
    }
  \end{equation}
  then there is an induced map $\phi^*: \V_1^\delta(E,n) \to \V_1^\delta(F,n)$ defined by
  $$
\phi^* ([x]_n) = [\phi^{(n)}(x)]_n.
  $$
In particular, if $E$ and $F$ are completely order isomorphic then $\V_1^\delta(E,n) \cong \V_1^\delta(F,n)$, for any $n \geq 1$.
\end{prop}

\begin{lma}
  \label{lma:direct-sum-V1}
Let $E$ and $F$ be unital operator system and let $\delta \geq  0$. Then for any $n$ we have $\V_1^\delta(E\oplus F,n) \cong \V_1^\delta(E,n) \times \V_1^\delta(F,n)$.
\end{lma}

Again, as with the group $K_0$ discussed in the previous section, we may introduce a semigroup structure on the direct limit $\V_1^\delta(E) := \varinjlim \V_1^\delta(E,n)$ of a direct system $(\V_1^\delta(E,n),\imath_{nm})$ with connecting morphisms defined as in \ref{eq:dir-syst}.
It is a semigroup when equipped with the direct sum $[x]+ [x'] = [x \oplus x']$ and identity element $0 = [e]$. 

\begin{defn}
  Let $(E,e)$ be a unital operator system and let $\delta \geq 0$. We define the $K$-theory group $K_1^\delta(E)$ of $E$ to be the Grothendieck group of $\V_1^\delta(E)$. 
\end{defn}
We will refer to this group as the {\em $\delta$-gapped $K_1$-group} of $E$. 
Similar to Proposition \ref{prop:stab} above we may show:
\begin{prop}[Stability of $K_1^\delta$]
  \label{prop:stab-K1}
  Let $E$ be a unital operator system, $\delta \geq 0$, and let $N$ be a natural number. Then $\V^\delta_1(E)$ is isomorphic to $\V^\delta_1(M_N(E))$ (and so are the corresponding $K_1^\delta$-groups).
\end{prop}

\section{Higher $\delta$-gapped K-groups and formal periodicity}
As in $K$-theory for $C^*$-algebras there are also higher-order invariants, which we now introduce. We first need the following.

\begin{defn}
  An operator system $E$ is called $\Z_2$-graded if there is a direct sum decomposition of $*$-vector spaces $E  = E^{(0)} \oplus E^{(1)}$ such that the corresponding grading $x \mapsto (-1)^{|x|}x $ is a complete order isomorphism.
\end{defn}
Consequently, there is a $*$-isomorphism on the $C^*$-envelope that extends the grading on $E$, turning $C^*_\env(E)$ into a $\Z_2$-graded $C^*$-algebra.

Much in the same way as Van Daele defined $K$-theory for graded Banach algebras \cite{Dae88a,Dae88b}, we will consider odd $\delta$-singular elements in $M_n(E)$ as our candidates for K-theory representatives. More precisely, we set

\begin{defn}
  Let $(E,e)$ be a $\Z_2$-graded unital operator system. 
  We define 
  $$
  \widehat H^{\delta}(E,n):= H^\delta(E,n) \cap M_n(E^{(1)}) ; \qquad (d=0,1).
  $$
and also set $\widehat \V^\delta(E,n) = \widehat H^{\delta} (E,n)/_{\sim_n}$.

\end{defn}

It would be interesting to develop this theory further for graded operator systems, similar to \cite{Dae88a,Dae88b}, however, we will focus here on the ungraded case. 
For this, we will exploit the graded Clifford algebras ({\em cf.} Appendix \ref{cliff}) to connect to our previous definitions of $\V_0^\delta$ and $\V_1^\delta$.
Let $E$ be an ungraded unital operator system, and consider the tensor product $E \otimes \CCl_p$. Given the structure of the Clifford algebras as matrix algebras, the operator system structure of this tensor product is unambiguous, and so is the grading on it.

\begin{defn}
  Let $(E,e)$ be an (ungraded) unital operator system and let $\delta \geq 0$. We define the higher invariants $\V_p^\delta(E,n)$ for $p=0,1,\ldots, $ by
  $$
  \V_p^\delta(E,n) :=\widehat \V_p^\delta (E \otimes \CCl_p,n) 
$$
The {\em higher $\delta$-gapped K-theory groups} $K_p^\delta(E)$ of a unital operator system $E$ are defined as the Grothendieck groups of the corresponding direct limits $\V_p^\delta(E) := \varinjlim \V_p^\delta(E, n)$, equipped with direct sum as the additive operation.
\end{defn}

  Let us now check that this definition is consistent with that of the invariants $\V_0^\delta (E,n)$ and $\V_1^\delta (E,n)$ (introduced in Definitions \ref{defn:class-herm} and \ref{defn:class-Gn}).

For the first case, $p=0$, we can indeed identify using Equation \eqref{eq:Cl2m-odd}:
  \begin{align*}
    H^\delta ( E \otimes \CCl_1^{(1)}  ,n ) &\cong \left \{ \begin{pmatrix} x & 0 \\ 0 & -x \end{pmatrix} \in H^\delta(E,2n) \right\}
    \cong H^\delta(E,n)
  . 
  \end{align*}
  The last isomorphism follows from the fact that there is a similarity transformation between the following matrices: 
  $$
\begin{pmatrix} s & 0 & x & 0 \\ 0 & s & 0 & -x \\ x & 0 & s & 0 \\ 0 & -x & 0 & s \end{pmatrix} \sim \begin{pmatrix} s & x \\ x & s \end{pmatrix} \otimes \mathbb I_2
  $$
  Upon taking the homotopy equivalence classes, we find agreement with Definition \ref{defn:class-herm}. Note that such a doubling has also been considered in the definition of K-theory of ungraded Banach algebras in \cite[Remark 2.13(iv)]{Dae88a}.

\bigskip
  
  For $p=1$ we have, using Equation \eqref{eq:Cl2m-even}:
 $$
  H^\delta ( E \otimes \CCl_2^{(1)},n ) \cong \left \{ \begin{pmatrix} 0 & x \\ x^* &0  \end{pmatrix} \in H^\delta(E,2n)
  \right\}
  \cong G^\delta(E,n),
  $$
  again via a similarity transformation, this time:
   $$
\begin{pmatrix} s & 0 & 0 & x \\ 0 & s & x^* & 0 \\ 0 & x & s & 0 \\ x^* & 0 & 0 & s \end{pmatrix} \sim \begin{pmatrix} s & x \\ x^* & s \end{pmatrix} \otimes \mathbb I_2
  $$
We thus find agreement of the above definition of $\V_1^\delta(E,n)$ with Definition \ref{defn:class-Gn}. 
  
More generally, we have the following, which may be considered as an operator system analogue of formal periodicity:
\begin{prop}
    \label{prop:formal-bott-V}
  Let $(E,e)$ be a unital operator system and let $\delta > 0$. Then 
  $$  \V_{2m}^\delta(E,n) \cong \V_0^\delta(M_{2^{m}}(E),n), \qquad \V_{2m+1}^\delta(E,n) \cong \V_1^\delta(M_{2^{m}}(E),n)$$
  for any $m \geq 0$.
  \end{prop}
\proof
For the even-dimensional case, we use \eqref{eq:Cl2m-odd} to obtain
\begin{align*}
  \V_{2m}^\delta(E,n)& =  H^{\delta} (E \otimes \CCl_{2m+1}^{(1)} ,n) /_{ \sim_n} \\
  & \cong  \left\{ \left(\begin{smallmatrix} x & 0 \\ 0 & -x \end{smallmatrix}\right) \in H^\delta(E \otimes( M_{2^m}(\C) \oplus  M_{2^m}(\C) ),n) \right\} /_{\sim_n}\\
  & \cong \left\{  x  \in H^\delta( M_{2^m}(E),n)\right\}/_{\sim_n} \equiv \V_0^\delta(M_{2^m}(E),n).
  %
  \intertext{For the odd-dimensional case, we instead use \eqref{eq:Cl2m-even} to write}
  \V_{2m+1}^\delta(E,n)& =  H^{\delta} (E \otimes \CCl_{2m+2}^{(1)},n )/_{\sim_n} \\
  & \cong \left\{ \left(\begin{smallmatrix} 0 & x \\ x^* & 0 \end{smallmatrix}\right) \in H^\delta(E \otimes M_{2^{m+1}}(\C),n)\right\}/_{\sim_n}\\
  & \cong \left\{  x  \in G^\delta(M_{2^m}(E),n) \right\}/_{\sim_n} \equiv \V_1^\delta(M_{2^m}(E),n).\qedhere
\end{align*}
\endproof
Proposition \ref{prop:formal-bott-V}, in combination with the stability of $K_0^\delta$ and $K_1^\delta$ (Proposition \ref{prop:stab} and \ref{prop:stab-K1}), then yields the following periodicity result:
\begin{thm}
Let $(E,e)$ be an ungraded unital operator system and let $\delta \geq 0$. Then $K_{2m}^\delta(E) \cong K_0^\delta(E)$ and $K_{2m+1}^\delta(E) \cong K_1^\delta(E)$.
  \end{thm}


\section{Application to the spectral localizer}
\label{sect:spectr-loc}
In \cite{LS18a,LS18b} the spectral localizer was introduced as a powerful tool for index pairings. They reduce the computation of a certain Fredholm index to the computation of the signature of a finite-dimensional matrix ---the spectral localizer. We will put it in the context of our notion of K-theory, realizing the spectral localizer as a spectral flow formula (very much as in \cite{LS18a,VSS19}) which maps $\V_p^\delta(E,n) \to \Z$. For more details on spectral flow, including the spectral localizer, we refer to the excellent textbook \cite{DSW23} and references therein.

In order to describe the index map from $\V_p^\delta(E,n)$, we need an operator system spectral triple \cite{CS20}. 
\begin{defn}
  A (unital) {\em operator system spectral triple} is given by a triple $(E,\H,D)$ where $E$ is a unital operator system realized concretely so that $E \subseteq C^*_\env(E) \subseteq B(\H)$, 
  and a self-adjoint operator $D: \dom(D) \to \H$ such that 
  \begin{itemize}
  \item the commutators $[D,x]$ extend to bounded operators for all $x \in \E$ for a dense $*$-subspace $\E \subseteq E$;
    \item the resolvent $(i+D)^{-1}$ is a compact operator.
  \end{itemize}
  An operator system spectral triple is called {\em even} if in addition to the above, there is a grading operator $\gamma$ (so that $\gamma^* = \gamma, \gamma^2 = 1_\H$) which commutes with all $x \in E$ and anti-commutes with $D$. Otherwise, it is called {\em odd}.
\end{defn}
In the even case, we can decompose $\H = \H_+ \oplus \H_-$ according to the eigenvalues of $\gamma$, and decompose accordingly
$$
D = \begin{pmatrix} 0 & D_0 \\ D_0^* & 0 \end{pmatrix}.
$$
In the odd case, there is no such grading, but for convenience we will then write $D_0 \equiv D$. 

Given an operator system spectral triple $(E,\H,D)$, a parameter $\kappa >0$ and an element $x \in G^\delta(\E,n)$ we now define the {\em generalized spectral localizer} for any $s \in (0,\delta)$ 
as the following (unbounded) self-adjoint operator on $\H^{\oplus 4n}$:
\begin{align}
  \L_{\kappa} (D,x,s)&= \begin{pmatrix} s & x &  \kappa D^{\oplus n}_0 & 0 \\ x^* & s & 0 & \kappa D_0^{\oplus n}  \\
 \kappa (D^{\oplus n}_0)^* & 0 &  - s & -x \\  0 & \kappa (D^{\oplus n}_0)^* & - x^* & -s  \\
  \end{pmatrix}.
\end{align}

We use the term 'generalized' to distinguish it from the even and odd spectral localizer consider in \cite{DSW23}. However, they are not unrelated:
\begin{rem}
  \label{rem:spectr-loc}
  Note that if we set $s=0$ in the odd case (so that $D_0 = D$) we obtain (twice) the odd spectral localizer $L_{\kappa}^{\text{odd}}$ introduced in \cite{LS18b}:
  $$
 \L_{\kappa} (D,x,s =0 ) \sim  (L^{\text{odd}}_\kappa)^{\oplus 2} ; \qquad L^{\text{odd}}_\kappa= \begin{pmatrix} \kappa D^{\oplus n} & x \\ x^* & - \kappa D^{\oplus n} \end{pmatrix},
 $$
 that is, up to a similarity transformation.
 
 In the even case, we will instead consider pairings with self-adjoint elements $x \in H^\delta(\E,n)$  in which case we find for $s=0$ the even spectral localizer $L^{\text{even}}_\kappa$ of \cite{LS18a}:
$$
\L_{\kappa} (D,x,s=0 ) \sim  (L^{\text{even}}_\kappa)^{\oplus 2} ; \qquad L^{\text{even}}_\kappa =\begin{pmatrix}  x & \kappa (D_0)^{\oplus n}  \\ \kappa (D_0^*)^{\oplus n} &  - x  \end{pmatrix}
$$
  \end{rem}

In general, the shift $s \in (0,\delta)$ is crucial in order to define the spectral localizer at the level of the operator systems, even though in some cases it may be set to zero ({\em cf.} Section \ref{sect:spectr-loc-S1} below). 
\begin{prop}
  \label{prop:homotopy-sf}
  Let $(E,\H,D)$ be a finite-dimensional operator system spectral triple. 
  Let $x \in G^\delta(\E,n)$ and let $s \in (0,\delta)$ and write $g_s = \min \{ s, \delta-s \}$. Then the
  signature of the generalized spectral localizer, 
     $$
  \Sig (\L_\kappa(D,x,s)),
       $$
  is constant in the $(\kappa,s)$-region defined by $0< \kappa < g_s^2 \| [D,x] \|^{-1}$ and $0< s < \delta$, and invariant under homotopy equivalence. Consequently, we have an induced map $\ind_D^\delta: \V_1^\delta(\E,n) \to \Z$ given by
  $$
\ind_D^\delta([x]) =
\lim_{\kappa,s \to 0} \frac 14  \Sig (\L_\kappa(D,x,s)),
$$
where the limit is taken in the above region.

     In the even case, we take $x \in H^\delta(\E,n)$ whence the induced map is $\ind_D^\delta: \V_0^\delta(\E,n) \to \Z$.
  \end{prop}

\proof
Without loss of generality we take $n=1$ and compute very similar to \cite{LS18a} that
\begin{align}
  \L_\kappa(D,x,s)^2 & = \begin{pmatrix}\left( \begin{smallmatrix}  \kappa^2 D_0 D_0^* & 0 \\0 &  \kappa^2 D_0 D_0^* \end{smallmatrix}\right) + \begin{pmatrix} s & x \\ x^* & s \end{pmatrix}^2
    &  \begin{pmatrix} 0 & - \kappa [D_0, x] \\ - \kappa[D_0,x^*] & 0 \end{pmatrix}\\
    \begin{pmatrix} 0 &  \kappa [D_0^*, x] \\  \kappa[D_0,x  ] & 0 \end{pmatrix}
    & \left( \begin{smallmatrix}  \kappa^2 D_0^* D_0 & 0 \\0 &  \kappa^2 D_0^* D_0 \end{smallmatrix}\right) + \begin{pmatrix} s & x \\ x^* & s \end{pmatrix}^2 \end{pmatrix}  \nonumber
  \\
  &\geq \left( g^2_s -\kappa \| [D,x]\| \right) 1_{M_4(B(\H))}.
\label{eq:loc-square}
\end{align}
Hence $\L_\kappa(D,x,s)$ is invertible (and thus has well-defined signature) provided $0 < \kappa< g_s^2 \| [D,x]\|^{-1}$ and $0< s < \delta$. Moreover, the signature is constant in this region since no eigenvalues will cross the origin. This also makes the limit of $\Sig (\L_\kappa(D,x,s))$ when $(\kappa,s)$ approach the origin from within this region well-defined.

Also, one may easily check that
$$
\Sig (\L_\kappa(D,x \oplus x',s)) = \Sig (\L_\kappa(D,x ,s))+ \Sig (\L_\kappa(D,x' ,s)).
$$
In order to see that $\Sig (\L_\kappa(D,e ,s)) = 0$ consider (in the odd case) an orthonormal eigensystem $\{v_\lambda\}_\lambda$ of $D$ in $\H$ so that $D v_\lambda = \lambda v_\lambda$. We may write the matrix of $\L_\kappa(D,x,s)$ in this eigensystem as
\begin{align}
 \langle v_\lambda, \L_\kappa(D,x,s) v_{\lambda'} \rangle  =& \delta_{\lambda\lambda'}
\begin{pmatrix} s & 1 & \kappa \lambda & 0\\
 1&  s & 0 & \kappa \lambda \\
 \kappa \lambda & 0 & -s & -1 \\
0 &  \kappa \lambda & -1 & -s \\
\end{pmatrix},
\label{eq:loc-basis}
\end{align}
which has set of eigenvalues $\{ \pm \sqrt{(1 \pm ' s)^2+\kappa^ 2 \lambda^2}\}$ so that the signature of the spectral localizer vanishes.

In the even case, we consider instead an eigensystem $\{v_\lambda =  (v_\lambda^+, v^-_\lambda)\}_\lambda$ in $\H_+ \oplus \H_-$ such that $D_0 v^-_\lambda = \lambda v^+_\lambda$ and $\gamma v_\lambda^\pm = \pm v_\lambda^\pm$. Again, we find that the matrix of $\L_\kappa(D,x,s)$ in this eigensystem is given by \eqref{eq:loc-basis} so that the signature of the spectral localizer vanishes. 

Consider now a homotopy $\tilde x$ in $G^\delta(E,n)$ (or in $H^\delta(E,n)$ in the even case) between $x = \tilde x(0)$ and $x' = \tilde x (1)$ with $s$-gap $\geq g_s$. By a computation similar to the one in Eq. \eqref{eq:loc-square} we find that 
$$
\L_{\kappa}(D,\tilde x,s)^2 \geq g^2_s - \kappa \sup_t \| [D, \tilde x(t)]\| 
$$
So as long as  $\kappa < g_s^2 / \sup_t \| [D,\tilde x(t)]\|$  we find that $\Sig \L_{\kappa}(D,\tilde x(t))$ is constant in $t$. We combine this with the fact that $\kappa <  g_s^2 / \| [D,\tilde x(0)]\|$ and $\kappa < g_s^2 / \| [D,\tilde x(1)]\|$ to conclude that $\lim_{\kappa, s \to 0} \Sig \L_{\kappa}(D,\tilde x(0),s) =\lim_{\kappa, s\to 0} \Sig \L_{\kappa}(D,\tilde x(1),s)$. This completes the proof. 
\endproof

We may now rephrase the main results of \cite{LS18a,LS18b}. For an invertible (self-adjoint) element $a \in A$ one considers the class $[a]$ in the K-theory group of the $C^*$-algebra $A$. Given a spectral triple $(A,\H,D)$, there is an index pairing $\ind_D([a])$. As shown in {\em loc.cit.} this index may be computed in terms of a spectrally truncated (operator system) spectral triple $(P A P,P \H , P D P)$ for a spectral projection $P$ (of $D$) of sufficiently high rank. Indeed, we then have
\begin{equation}
  \ind_D ([a]) = \frac 14 \Sig(\L_\kappa (P DP ,x, 0)) = \frac 12 \Sig(L_\kappa^{\text{even/odd}}) 
 \label{eq:ind-loc}
\end{equation}
since in fact it is shown in \cite{LS18b} that $\L_\kappa (D,x, 0)$ already has a positive gap, say $\L_\kappa (D,x, 0)^2 \geq  g^2$. Thus, adding the matrix with diagonal entries $(s,s,-s,-s)$ to obtain $\L_\kappa (D,x, s)$ will not change this property, as long as $s<g$.

\subsection{Example: spectral localizer on the circle}
\label{sect:spectr-loc-S1}
Let us illustrate the spectral localizer by spectral truncations of the circle. We consider a unitary $u(t) = e^{imt}$ with winding number $m$ and want to use the spectral localizer to compute the index pairing with $D_{\mathbb S^1} = -i d/dt$.

We take a spectral projection $P_N$ onto $\text{span}_\C \{ e_{-N}, e_{-N+1}, \ldots, e_N \}$ and realize that $P_N u P_N$ is a Toeplitz matrix with $1$'s on the $m'th$ diagonal, and zeros elsewhere. Note that $P_N u P_N$ is not invertible, however the following self-adjoint $2 \times 2$ matrix with values in the truncated operator system {\em is} non-singular for any $s \in (0,\delta)$ for some $\delta>0$:
$$
\begin{pmatrix} s P_N &  P_N u P_N \\ P_N u^* P_N &  s P_N \end{pmatrix} 
$$
Hence $P_N u P_N  \in G^\delta(P_N C(\mathbb S^1)P_N,1)$.

\begin{figure}
  \includegraphics[scale=.38]{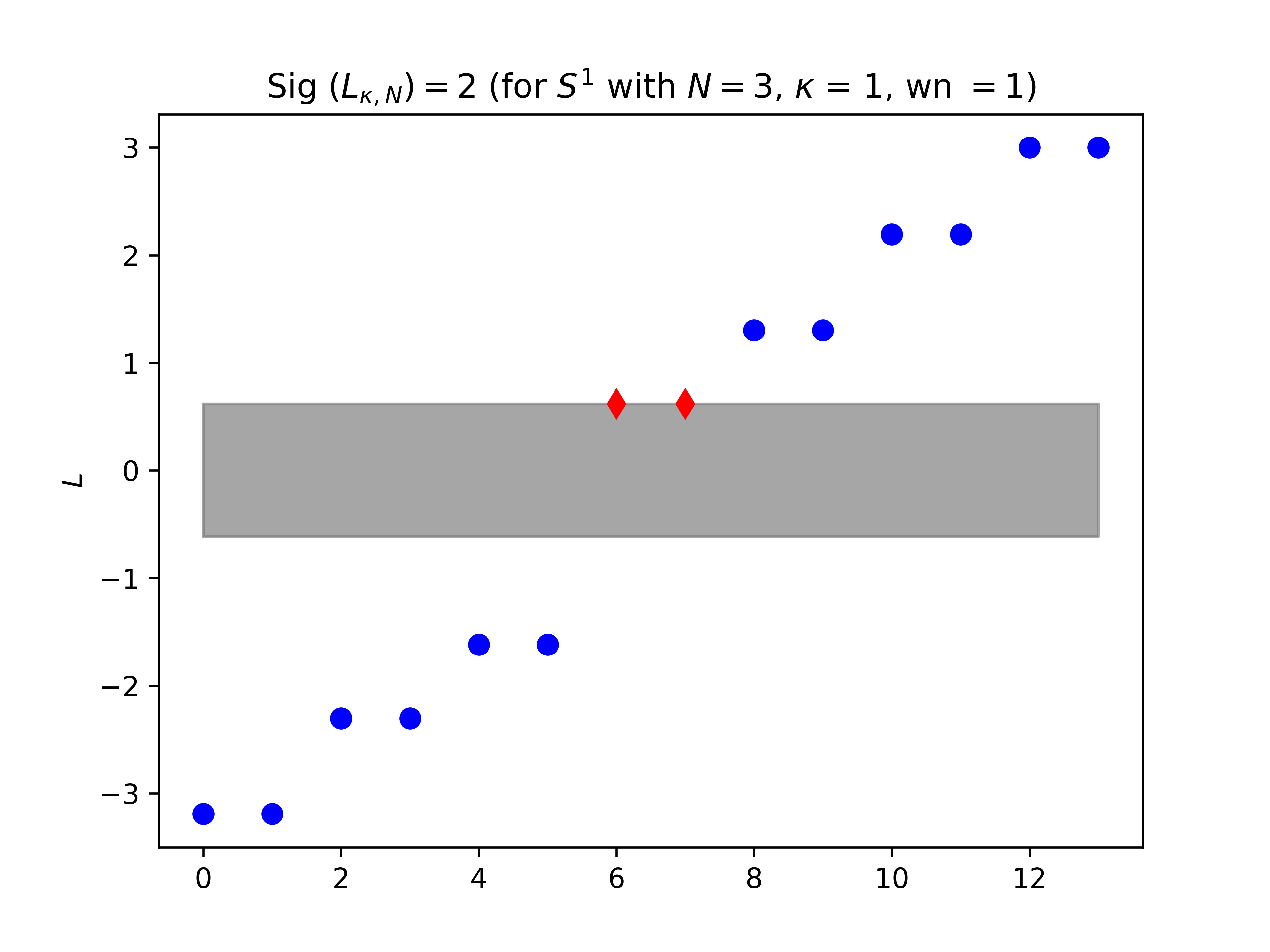}
  \includegraphics[scale=.38]{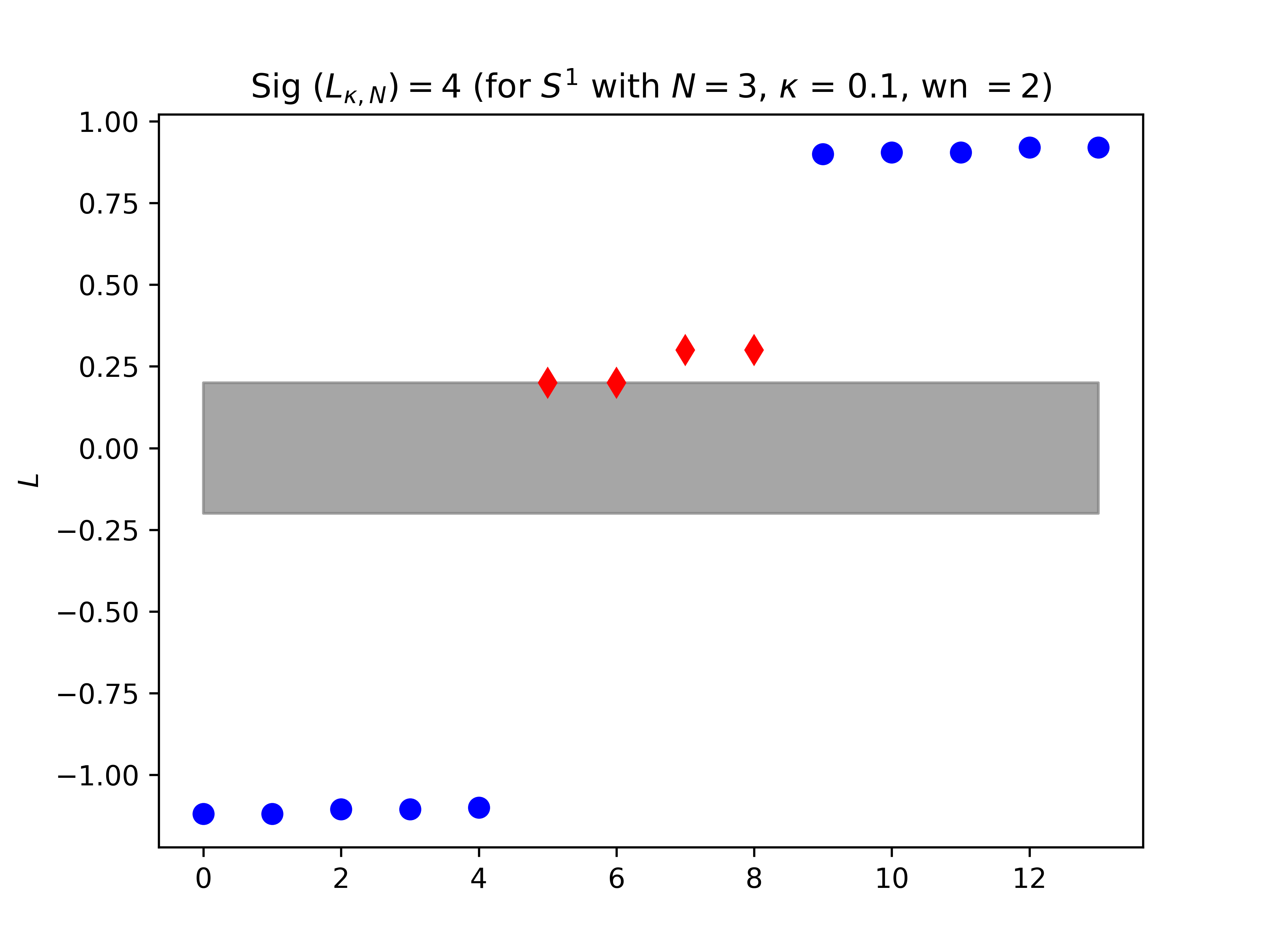}
  \caption{The negative and positive eigenvalues for the spectral localizer $L_{\kappa,N} = L_{\kappa,N,s=0}$ on the circle for $m=1, \kappa= 1, N=3$ (left) and  $m=2, \kappa= 0.1, N=3$ (right). The red (diamond-shaped) dots indicate the surplus of positive eigenvalues as compared to the negative ones.}
    \label{fig:L-S1}
\end{figure}

The spectral localizer is given in this case by the following matrix: 
  $$
\L_{\kappa,N,s}:= \L_{\kappa}(P_N DP_N, P_N u P_N, s) = \left(\begin{smallmatrix} s & P_N u P_N &   \kappa P_NDP_N & 0 \\  P _N u^* P_N & s & 0 &  \kappa P_NDP_N  \\
  \kappa P_NDP_N & 0 &  - s & -P_N u P_N  \\  0 &  \kappa P_NDP_N & -  P _N u^* P_N & -s  \\
\end{smallmatrix} \right).
$$
The main result of \cite{LS18b} ({\em cf.} Eq. \eqref{eq:ind-loc} and the discussion following it) can now be used to show that for suitable $N,\kappa$ there exists $g>0$ so that the generalized spectral localizer satisfies $\L_{\kappa ,N,s}^2 \geq g^2$ for any $s< \delta$, including $s =0$. As such, the signature can be computed by setting $s =0$, so that 
$$
\ind_{P_N DP_N }([P_N u P_N ]) =
\frac 14 \text{Sig} ~\L_{\kappa ,N,{s=0}} = \frac 12 \text{Sig} ~ L_{\kappa,N},
$$
where as in Remark \ref{rem:spectr-loc} we have written $\L_{k,N,0} = L_{\kappa,N} \oplus L_{\kappa,N}$ in terms of
$$
L_{\kappa, N} = 
\begin{pmatrix} \kappa P_N D P_N &  P_N u P_N \\ P_N u^* P_N &   -P_N D P_N \end{pmatrix} $$
The main result of \cite{LS18b} states that $\ind (P_D u P_D) = \frac 12\text{Sig} ~  L_{\kappa,N}$ so that we deduce that
$$
\ind (P_D u P_D) = \ind_{P_N DP_N }([P_N u P_N ]) 
$$

In Figure \ref{fig:L-S1} we illustrate the resulting signature by showing the negative and positive eigenvalues of $L_{\kappa,N}$. We find already for low $N$ that the spectral localizer is equal to (twice) the winding number of $u$.

\appendix

\section{Clifford algebras}
\label{cliff}
Recall the definition of the Clifford algebras $\CCl_p$ for $p \geq 1$ associated to the vector spaces $\C^p$  with their standard quadratic form. 
These are $\Z_2$-graded $C^*$-algebras, {\em i.e.} $\CCl_p = \CCl_p^{(0)} \oplus \CCl_p^{(1)}$. For each $p$ we have a grading operator $\Gamma$ such that $\CCl_p^{(d)} = \{ a \in \CCl_p: \Gamma a = (-1)^d a \Gamma \}$ (if $p$ is odd, this grading operator is an element in the multiplier algebra of $\CCl_p$). 

In the odd case $p=2m+1$ we have $\CCl_{2m+1} \cong M_{2^m}(\C) \oplus M_{2^m}(\C)$ and we may take
$$
\Gamma = \begin{pmatrix} 0 &  \mathbb I_{2^{m}} \\  \mathbb I_{2^{m}}  & 0 \end{pmatrix} \in M_{2^{m+1}}(\C).
$$
Consequently, 
 \begin{align}
  \CCl_{2m+1}^{(0)}&\cong\left\{ (x,x) :  x \in M_{2^{m}}(\C) \right\} \subset M_{2^m}(\C) \oplus M_{2^m}(\C).\nonumber \\
\CCl_{2m+1}^{(1)} & \cong\left\{ (x,-x) : x \in M_{2^{m}}(\C) \right\}\subset M_{2^m}(\C) \oplus M_{2^m}(\C) .\label{eq:Cl2m-odd}
\end{align}

 Instead, if $p=2m$ then $\CCl_{2m} \cong M_{2^m}(\C)$ and we may take
$$
\Gamma = \begin{pmatrix} \mathbb I_{2^{m-1}} & 0 \\ 0 & - \mathbb I_{2^{m-1}} \end{pmatrix}.
$$
Accordingly, 
\begin{align}
  \CCl_{2m}^{(0)}& \cong\left\{ \begin{pmatrix} x & 0 \\ 0 & y \end{pmatrix}:  x,y \in M_{2^{m-1}}(\C) \right\}. \nonumber \\
\CCl_{2m}^{(1)} &\cong\left\{ \begin{pmatrix} 0 & x \\ y & 0 \end{pmatrix}:  x,y \in M_{2^{m-1}}(\C) \right\}.\label{eq:Cl2m-even}
\end{align}

\newcommand{\noopsort}[1]{}\def\cprime{$'$}


\end{document}